\documentclass[reqno, 12pt]{amsart}

\usepackage[utf8]{inputenc}

\usepackage{a4wide}
\usepackage{color}
\usepackage{amsmath,amsfonts,latexsym,url}

\usepackage{bm}
\usepackage[colorlinks=true,citecolor=blue]{hyperref}
\usepackage{xcolor}

\usepackage{todonotes}

\newtheorem{definition}{Definition}

\newtheorem{theorem}[definition]{Theorem}

\makeatletter
\def\iddots{\mathinner{\mkern1mu\raise\p@
\vbox{\kern7\p@\hbox{.}}\mkern2mu
\raise4\p@\hbox{.}\mkern2mu\raise7\p@\hbox{.}\mkern1mu}}
\makeatother

\def\K{{\mathbb K}}

\def\Qbar{\overline{\mathbb{Q}}}
\def\sing{\textup{Sing}}
\def\trs{\mathcal{S}}

\makeatletter
\newcommand*{\house}[1]{%
 \mathord{%
 \mathpalette\@house{#1}%
 }%
}
\newcommand*{\@house}[2]{%
 \dimen@=\fontdimen8 %
 \ifx#1\scriptscriptstyle\scriptscriptfont
 \else\ifx#1\scriptstyle\scriptfont
 \else\textfont\fi\fi
 3 %
 \sbox0{%
 $#1%
 \vrule width\dimen@\relax
 \overline{%
 \kern2\dimen@
 \begingroup 
 #2%
 \endgroup
 \kern2\dimen@
 }%
 \vrule width\dimen@\relax
 \mathsurround=1.5\dimen@ 
 $%
 }%
 \ht0=\dimexpr\ht0-\dimen@\relax
 \dp0=\dimexpr\dp0+2\dimen@\relax
 \vbox{%
 \kern\dimen@ 
 \copy0 %
 }%
}

\begin{document}

\title[]
{Explicit degree bounds for right factors \\ of linear 
differential operators}

\author[]{A.~Bostan}
\address{Inria, Universit\'e Paris-Saclay, 1 rue Honor\'e d'Estienne d'Orves, 91120 Palaiseau, France}

\author[]{T.~Rivoal}
\address{Institut Fourier, CNRS et Universit\'e Grenoble Alpes, CS 40700, 38058 Grenoble cedex~9} 
\author[]{B.~Salvy}
\address{Univ Lyon, EnsL, UCBL, CNRS, Inria,  LIP, F-69342, LYON Cedex 07, France}

\begin{abstract} 
If a linear differential operator with rational function coefficients is
reducible, its factors may have coefficients with numerators and denominators
of very high degrees. When the base field is $\mathbb C$, we give a completely explicit bound for the degrees of
the monic right factors in terms of the degree and the order of the original
operator, as well as of the largest modulus of the local exponents at all its
singularities. As a consequence, if a differential operator~$L$ 
has rational function coefficients over a number field, we get degree bounds for its monic right factors in terms of the degree, the order
and the height of~$L$, and of the degree of the number field.

\end{abstract}

\date{\today}



\maketitle

\date{\today}

\section{Introduction}

\paragraph{\em Context}

We are interested in factorizations of linear differential operators in
$\mathbb K(z)[\frac{d}{dz}]$, where $\mathbb K$ is either $\mathbb C$ or
$\Qbar$ (embedded into $\mathbb C$). In the latter case, there is no loss of generality in assuming that
$\mathbb K$ is a number field (because the coefficients all live in such a
number field) and in this case we denote its degree by $\kappa:=[\mathbb K :
\mathbb Q]$.

Without loss of generality, we assume that $L\in \mathbb K[z][\frac{d}{dz}]$, i.e., it has the form
\begin{equation}\label{eq:1}
L=\sum_{j=0}^{m} p_j(z)\Big(\frac{d}{dz}\Big)^j
\end{equation}
for some polynomials $p_j(z)\in \mathbb K[z]$, with
$p_m\neq0$. The order of $L$ is $m$, and we assume that $m\ge 1$.

Assume that there exists a factorization $L=NM$ with $M,N \in \mathbb K(z)[\frac{d}{dz}]$, where 
\[
M=\sum_{j=0}^{r} \frac{A_j(z)}{B(z)}\Big(\frac{d}{dz}\Big)^j,
\]
for polynomials $A_0(z), \ldots, A_r(z)$ and $B(z)$ in $\mathbb K[z]$, 
with $A_r\neq0$ and $B$ of minimal degree. The order of $M$ is $r$. 
We call \emph{degree} of $M$ the quantity
\[\deg_z(M):=\max(\deg(A_0),\dots,\deg(A_r),\deg(B)).\]

Obviously $r\le m$, because the order of $L$ is the sum of the orders of $N$ and $M$.
But it is well known that $\deg_z(M)$ can be much larger than~$q:=\deg_z(L)$,
and it is in fact notoriously difficult to control~$\deg_z(M)$ in terms of $L$. To the
best of our knowledge, the first (and so far the only)  
bound for~$\deg_z(M)$ has been given by Grigoriev \cite[Theorem~1.2]{Grigoriev90}. 
On the one hand, Grigoriev's bound holds for any factor, not only
for right factors. But on the other hand, it is only an asymptotic bound; for
instance, with respect to the input degree~$m$, it writes $\exp \left( 2^{m
\cdot o(2^m)} \right)$ as $m\to +\infty$. One of our aims is to replace this bound by an \emph{effective bound}, i.e., without any constant implicit in a $o()$- or $O()$-estimate. This is important to ensure the termination of a recent algorithm by Adamczewski and the second author~\cite{AdRi18}; see below.
\medskip

\paragraph{\em Main result}
Here, we seek entirely explicit bounds holding for all~$m$ and for {\em any}
operator $L \in \Qbar[z][\frac{d}{dz}]$. As we will see, such a bound is a
consequence of the following result.

\begin{theorem} \label{theo:1} 
Let $L\in \mathbb C[z][\frac{d}{dz}]$ and $M$ be a monic right factor of $L$. Then the degree of $M$ satisfies
\begin{equation}\label{eq:n16011}
\deg_z(M)\le 
r^2(\trs+1)\mathcal{E}+ r(\mathcal{N}+1)\trs + r\mathcal{N} +\frac12 r^2(r-1)\big((\trs+1)(\mathcal{N}+1)-2\big),
\end{equation}
where 
\begin{itemize}
\item[$\bullet$] $r \ge 0$ is the order of~$M$;
\item[$\bullet$] $\mathcal{E}\ge 0$ is the largest modulus of the local generalized exponents of~$L$ at $\infty$ and at its finite non-apparent singularities; 
\item[$\bullet$] $\mathcal{N}\ge 0$ is the largest of all the slopes  of $L$ at its finite singularities and at $\infty$; 
\item[$\bullet$] $\trs \ge 0$ is the number of finite non-apparent  singularities of $L$.
\end{itemize}
\end{theorem}

The notions of apparent singularities, generalized exponents and slopes of a differential operator are
recalled in \S\ref{sec:theo1}. 
$L$ is Fuchsian if and only if $\mathcal{N}=0$, in which case the generalized exponents are the usual exponents of regular or regular singular points. 
If $\mathcal{N}\ge 1$, then on the right-hand side of
\eqref{eq:n16011}, the term $r\mathcal{N}$ can be replaced by $r(\mathcal{N}-1)$; see the discussion
following inequality~\eqref{eq:t100} in \S\ref{ssec:casgeneral}. Also, the factor $(\mathcal{S}+1)\mathcal{E}$ in the first term of the right-hand side of~\eqref{eq:n16011} can be further refined and replaced by the sum of the largest moduli of the local generalized exponents of~$L$ at~$\infty$ and at all its finite non-apparent singularities, rather than $(\mathcal{S}+1)$ times their maximum value.
We refer the reader to the comments made after the proof of Theorem~\ref{theo:1} concerning the choice of $\mathbb C$ as the base field instead of an arbitrary algebraically closed field of characteristic~0.

\medskip

\paragraph{\em Bounding the degree of $M$ in terms of the degree~$q$, the order~$m$
and the height~$H$ of~$L$}

Note that $r\le m$, $\mathcal{N}\le m+q$ and $\trs \le q$, so that Theorem~\ref{theo:1} reduces the
problem of bounding~$\deg_z(M)$ to the determination of an explicit upper bound for
$\mathcal{E}$, or rather for the larger quantity~$E$ defined as the largest modulus of the local generalized exponents of~$L$ at~$\infty$ and at its finite  singularities. Bounds for $E$ are known in the case where $L\in
\K(z)[\frac{d}{dz}]$, where $\K$ is a number field of degree~$\kappa$, embedded into $\mathbb C$. 
Grigoriev exhibited such a bound in that case, but again his
result~\cite[Corollary, p.~21]{Grigoriev90} is only asymptotic in the
order~$m$ of~$L$, see below. The first entirely explicit bound for
$E$ was obtained in 2004 by Bertrand, Chirskii and
Yebbou~\cite{BeChYe04}. Their approach was based in part on Malgrange's
truncation method, which was eventually published in
\cite[pp.~97--107]{DeMaRa}. In terms of the height $H$ of
the operator ~$L$~\cite[p. 246 and p. 252]{BeChYe04}, their bound reads
\begin{equation}\label{eq:19113} 
E\le
2^{(36(q+1)m\kappa)^{9(q+1)^2 m^{3m}}}H^{(5\kappa(q+1)m)^{9(q+1)^2
m^{3m}}}. 
\end{equation}

The inequalities~\eqref{eq:n16011} and~\eqref{eq:19113}, together with the
bounds $r\le m$, $\mathcal{N}\le m+q$ and $\trs \le q$, completely solve the problem of
finding an explicit upper bound for the degree of any monic right factor $M$
of~$L$, when $L\in \Qbar(z)[\frac{d}{dz}].$ It seems to be the first of this
type in the literature. We have chosen to formulate Theorem \ref{theo:1} in
terms of $\mathcal{E}$ as a parameter because the upper bound in
\eqref{eq:19113} seems pessimistic and any improvement of it would implicitly
improve Theorem \ref{theo:1}. On the other hand, the other terms on the
right-hand side of \eqref{eq:n16011} are already polynomial in the parameters
and are thus probably only marginally improvable.

\medskip

\paragraph{\em Asymptotic comparisons}

Below, we let $\mathcal{P}(X)$ denote different polynomials in
$\mathbb{Z}[X]$, with degree and coefficients independent of $\kappa,
m$ and $q$. With our notations, Grigoriev obtains the asymptotic estimate
$E\le H^{\mathcal{P}(\kappa q m)^m}$ as $m\to +\infty$, which is
much better than~\eqref{eq:19113}, which reads
$E\le H^{\mathcal{P}(\kappa q m)^{q^2 m^{3m}}}$ as $m\to
+\infty$. When $L$ is Fuchsian, Grigoriev's method as well as that of Bertrand \emph{et al.}
\cite[p. 254]{BeChYe04} provide better bounds, which turn out to be both of the form
$E\le H^{\mathcal{P}(\kappa q m)}$; one may wonder if this is asymptotically
optimal as $m\to +\infty$. In the general case, it would obviously be
interesting to close the gap between the uniform bound \eqref{eq:19113} and
Grigoriev's asymptotic bound for $E$. It would also be
interesting to do so in intermediate cases where some properties of $L$ are known
in advance. For instance, for applications related to $E$-functions (see
\cite{AdRi18}), $L$ may have only two singularities: $z=0$ which is regular,
and $z=\infty$ which is irregular with slopes in $\{0,1\}$.

\medskip

\paragraph{\em Optimality of the bound in Theorem~\ref{theo:1}} For any
integer $k\geq 1$, the second-order operator $L := z \left( \frac{d}{dz}
\right)^{2} + (2-z) \frac{d}{dz} +k$ admits the right factor $M = \frac{d}{dz}
- \frac{H'(z)}{H(z)}$, where $H(z)$ is the confluent hypergeometric Kummer
polynomial $H(z) = {_1 F _1} (-k;\,2;z) = \sum_{\ell=0}^k \binom{k}{\ell}
\frac{(-z)^\ell}{(\ell+1)!} $. Thus, $m=2$, $q=1$, $r=1$ and $\deg_z(M)=k$, and it is
easy to check that $\mathcal{E}=k$, $\mathcal{N}=1$ and $\trs=0$. Therefore the
bound of Theorem~\ref{theo:1} writes $\deg_z(M) \leq k$ (using the above mentioned improvement in the case $\mathcal{N}\ge 1$). The 
bound~\eqref{eq:n16011} is thus optimal for this example.

\medskip

\paragraph{\em Degrees of left factors.} Taking formal adjoints exchanges left
and right factors: if $L = NM$, then $L^\star = M^\star N^\star$, see
e.g.~\cite[p.~39--40]{putsinger}. Therefore, one can effectively bound the
degree of the left factor $N$ as well, by applying Theorem~\ref{theo:1} to
$N^\star$ and using the fact that all the quantities (order, degree, largest
slope, maximal exponent modulus, number of finite non-apparent singularities), involved in the inequality~\eqref{eq:n16011}
for $L^\star$ and $N^\star$ can be expressed or bounded in terms of the same set of quantities for $L$ and~$N$.

\medskip

\paragraph{\em Minimal differential equations} 
Besides its own interest, one of our motivations to study this factorization
problem comes from combinatorics~\cite{BoBoKaMe16} and number theory \cite{AdRi18,FiRi19}, where certain
D-finite power series in $\Qbar[[z]]$, called $E$- and $G$-functions, are
under study. In both cases, it is useful to be able to
perform the following task efficiently: given $f(z)\in\Qbar[[z]]\setminus \{0\}$ and $L\in
\Qbar(z)[\frac{d}{dz}]\setminus \{0\}$ such that $Lf(z)=0$, determine $M\in
\Qbar(z)[\frac{d}{dz}]\setminus \{0\}$ such that $Mf(z)=0$ and $M$ is of
minimal order with this property. Obviously, $M$ is then a right factor of $L$
and Theorem \ref{theo:1} applies to it. Assume $L\in \mathbb
K(z)[\frac{d}{dz}]$ with the same data as above and $\mathbb K$ a number
field, and let $f(z)\in \mathbb K[[z]]$ be a solution of the differential
equation $Ly(z)=0$. The power series $f$ need not be convergent. For any
integers $r, n$ such that $1\le r \le m$ and $n\ge 0$, define
$$
R(z):= \sum_{j=0}^r P_j(z)f^{(j)}(z),
$$
where $P_j(z)\in \mathbb K[z]$ are all of degree at most $n$. Then
$R(z)=\sum_{k=0}^\infty r_k z^k$ is a formal power series in $\mathbb K[[z]]$,
and we denote by $N$ its valuation (or order) at $z=0$, i.e., $N$ is the
smallest integer $k\ge 0$ such that $r_k\neq 0$. 
A key inequality is the following upper bound on~$N$~\cite{BeChYe04}: either  
$R(z)$ is identically zero or
\begin{equation}\label{eq:19111}
N\le r(n+1) + 2(q+1)^2m^3 + 2(q+1)m^2(E+1).
\end{equation}
This is proved by putting together results by Shidlovskii~\cite[Lemma~8, p.~83 and Eq.~(83), p.~99]{shidlovskii} and Bertrand, Chirskii, Yebbou~\cite[Thm.~1.2, p.~245]{BeChYe04}. With our notations, this yields
$N\le r(n+1)+n_0$ where $n_0$ is a quantity bounded above by $2(q+1)m^2(\mathcal{R}+1),$ with $\mathcal{R}\le E+(q+1)m$, see \cite[p.~252]{BeChYe04}.

Now, given $n$ and $r+1$ polynomials $P_j$, not all zero, letting $N$ denote the upper bound in Eq.~\eqref{eq:19111}, if the first $N+1$ Taylor coefficients of $R(z)$ are all 0, then $R(z)$ is proven identically zero, which means that $f(z)$ is a solution of 
$$
M:=\sum_{j=0}^r P_j(z) \Big(\frac{d}{dz}\Big)^j \in \mathbb K(z)\Big[\frac{d}{dz}\Big] \setminus\{0\},
$$
and thus $M$ is a right factor of~$L$. 

This remark was used by Adamczewski and the second author \cite{AdRi18}
to give an algorithm that
computes a non-zero operator $M$ such that $Mf(z)=0$ and $M$ is of minimal
order with this property. The input is $L\in \mathbb K(z)[\frac{d}{dz}]$ and
sufficiently many initial Taylor coefficients of $f$, so that the following
ones can be computed using $L$. Let $\widehat{n}$ be the quantity on the right-hand
side of the inequality~\eqref{eq:n16011}. The algorithm first sets $r=1$ and looks for $R$ with order $r$
and degree $\lceil \widehat{n} \rceil$ by requiring that its first $N+1$ Taylor
coefficients all be $0$ (this amounts to solving a homogeneous linear system
with algebraic coefficients given by the Taylor coefficients of~$f$). If no
non-zero solution is found, $r$ is increased and the same procedure is repeated, and
so on up to $r=m$ if necessary. In the end, $M\neq 0$ minimal for $f$ will be
found.

This algorithm  is not very efficient in
practice. Moreover, the inequalities~\eqref{eq:19113} and~\eqref{eq:19111}, as well as Grigoriev's Theorem~1.2 are all used to ensure the termination of the algorithm. It is important however to use Theorem~\ref{theo:1}  instead of Grigoriev's, as it holds for
arbitrary differential operators $L$ and $M \in \Qbar(z)[\frac{d}{dz}]$,
and not only asymptotically. 
Rather than using a uniform \emph{a priori} bound, a much more efficient minimization
algorithm, computing tight bounds dynamically along the lines of this article and van~Hoeij's work~\cite[Sec.~9]{Hoeij97b}, is under development~\cite{BoRiSa2020}.

\medskip

\paragraph{\em Related works.} The proof of Theorem \ref{theo:1} does not use
Grigoriev's method~\cite{Grigoriev90}, which relies on a subtle analysis of
Beke's classical factorization algorithm \cite[p. 118, \S 4.2.1]{putsinger},
see also~\cite{Singer81}. Instead, our method is inspired by van Hoeij's
factorization algorithm~\cite{Hoeij97, Hoeij97b}. This algorithm  
internally computes, on any input operator~$L$, 
upper bounds for the number of apparent singularities and for the degree of right factors of~$L$, using the
generalized Fuchs relation between local exponents. 
He did 
not give any explicit {\em a priori} degree bound, valid for any operator~$L$. Our main contribution here is such a bound when  
the base field is a number
field. It is difficult to trace back exactly when in the 80's the
(generalized) Fuchs relation was found to be relevant in this type of problem. In the Fuchsian case, it was used by Chudnovsky \cite{Chudnovsky80} to bound the number of apparent singularities in order  
to obtain an effective multiplicity estimate. See also \cite[p. 364, Example 2.7]{cormier} for a similar use of Fuchs' relation. 
Chudnovsky's multiplicity estimate was adapted by 
Bertrand and Beukers to the general case with the help of the generalized
Fuchs relation~\cite{BeBe85}. They obtained a multiplicity estimate in which
the effectivity of one specific constant was not completely clear. This
effectivity issue was eventually solved by Bertrand, Chirskii and Yebbou
\cite{BeChYe04}.

\section{Proof of Theorem \ref{theo:1}} \label{sec:theo1}

From this point on, we write $\partial_z$ for $\frac{d}{dz}$. Consider a monic operator
\begin{equation}\label{eq:monicR}
R=\sum_{j=0}^\mu c_j(z) \partial_z^j=\sum_{j=0}^\mu\frac{U_j(z)}{V(z)}\partial_z^j \in \mathbb C(z)[\partial_z],
\end{equation} 
with $U_j,V\in \mathbb C[z]$ and $V$ of lowest degree. We have $U_\mu=V$, $\deg(c_j):=\deg(U_j)-\deg(V)$ and $\deg_z(R):=\max(\deg(U_0), \ldots, \deg(U_{\mu-1}), \deg(V))$. By definition, the set $\sing(R)$ of 
finite singularities of~$R$ is the set of roots of $V$.  
Amongst the finite singularities of
$R$, we denote by $\alpha(R)$ the set of the \emph{apparent} ones, i.e., those
at which $R$ admits a local basis of power series solutions. Note that an apparent singularity is necessarily regular.
We denote by
$\sigma(R)$ the set of finite singularities of $R$ which are not in
$\alpha(R)$, so that $\sigma(R)$ and $\alpha(R)$ form a partition of
$\sing(R)$. For an operator with polynomial coefficients such as~$L$ from~\eqref{eq:1}, the sets $\sing(L)$, $\alpha(L)$ and $\sigma(L)$ are defined as the corresponding sets  for $R:=(1/p_m)\cdot L$.
In a factorization $L=NM$, we have $\sigma(M) \subset \sigma(L)
\subset \textup{Sing}(L)$ but $\alpha(M)$ may have no common element with
$\textup{Sing}(L)$. Because of this, the main difficulty in the method
presented below is to bound the number of apparent singularities of a
right factor of $L$.

We split the proof of the theorem into two parts. We start with the Fuchsian case because it is simpler but at the same time it  contains essentially all the ideas needed to prove the general case. We refer to the book by van der Put and Singer~\cite[\S4.4]{putsinger} for the definitions of classical notions related to linear differential operators (regular singularity, Fuchsian operator, local exponents,\dots).

\subsection{Fuchsian case}

Assume that we have a factorization $L=NM$ with operators $N,M$ in
$\mathbb C(z)[\partial_z]$ for which the operator~$M$ is Fuchsian and monic.
Note that $L$ need not necessarily be Fuchsian itself. We compute an explicit
upper bound on $\deg_z(M)$ in terms of $\mathcal{E}$. Our strategy is
inspired by van~Hoeij's approach~\cite{Hoeij97b}, itself based on ideas by
Chudnovsky~\cite{Chudnovsky80} and Bertrand-Beukers~\cite{BeBe85}, see
also~\cite{Ohtsuki82}.

The Fuchsianity of $M$ implies that it can be written
\[M = \partial_z^r + \frac{A_{1}(z)}{A(z)}\partial_z^{r-1} + \cdots + \frac{A_r(z)}{A(z)^r}, 
\]
where $A(z)$ has simple roots and $\deg(A_{i}) \leq \deg(A^i)-i$; see~\cite[Chap.~V, \S20, p.~77]{Poole1960}. All we now have to do is to derive an upper bound on the degree of $A$. This is done in two steps. The polynomial $A(z)$ can be factored in $\mathbb{C}[z]$ as $A(z) =
A_{\textsf{sing}}(z)A_{\textsf{app}}(z)$, where the roots of
$A_{\textsf{sing}}$ are the elements of $\sigma(M)$, while those 
of  $A_{\textsf{app}}$ are the elements of $\alpha(M)$.
Since $A_{\textsf{sing}}$ has simple roots, its degree is equal to~$\#\sigma(M)\le \#\sigma(L)=\trs$. 

The degree of $A_{\textsf{app}}(z)$ is equal to $\#\alpha(M)$ and it can be bounded above using  Fuchs' relation, which we now recall. We set 
\begin{equation}\label{eq:fuchs2}
S_\rho({M}) := \sum_{j=1}^r e_j(\rho) - \frac{r(r-1)}{2}
\end{equation}
where the $e_j(\rho)$'s are the local exponents of~$M$ at the
point~$\rho$, so that clearly $S_\rho({M})=0$ 
 when $\rho \in \mathbb C \cup \{\infty\}$ is an ordinary point of~$M$.  Then Fuchs' relation is the following global relation (see \cite[p.~138]{putsinger} for a different but equivalent version):
\begin{equation}\label{eq:Fuchs}
\sum_{\rho \in \mathbb{C} \cup \{\infty\}} S_\rho({M}) = \sum_{\rho
\in \textup{Sing}(M)\cup \{\infty\}} S_\rho({M}) = \; - r(r-1).
\end{equation}	

Now, the main observation is that if~$\rho\in \alpha(M)$, then $S_\rho({M}) \in
\mathbb{N} \setminus \{ 0 \}$  \cite[Chap.~V, \S18,
p.~69]{Poole1960}
\footnote{{\em Stricto sensu}, \cite{Poole1960} proves this under an {\em a priori}  
stronger definition of an apparent singularity~$\rho$, which requires the
\emph{holomorphy} of the basis of solutions at~$\rho$. Note, however, that the
proof is algebraic and does not use this assumption, see also~\cite[p.~187--188]{putsinger}.},
so that
$$
\#\alpha(M) \le \sum_{\rho \in \alpha(M)} S_\rho(M)
$$
and by~\eqref{eq:Fuchs} this implies that 
\begin{equation} \label{eq:newreferee}
\#\alpha(M) \leq -r (r-1) - \sum_{\rho \in \sigma(M)\cup \{\infty\}
}S_\rho(M).
\end{equation}

Since $M$ is a right divisor of~$L$, we have $\sigma(M) \subset \sigma(L)
$  and for any such singularity~$\rho \in \sigma(M)$, the exponents of $M$
at~$\rho$ are also exponents of~$L$ at~$\rho$, so that  
$|S_\rho({M})| \leq r \mathcal{E} + r(r-1)/2$ by \eqref{eq:fuchs2}.
Therefore, 
\[ \#\alpha(M) + r (r-1) \leq \sum_{\rho 
\in \sigma(L)\cup \{\infty\}} |S_\rho(M)| \leq (\trs+1)\left(r \mathcal{E} + \frac{r(r-1)}2 \right)\]
and 
\[ \#\alpha(M) \leq r(\trs+1) \mathcal{E} + \frac12r(r-1) (\trs-1).\]
Hence, 
\[ \deg(A) \leq r(\trs+1)  \mathcal{E} + \frac12 r(r-1) (\trs-1) +\trs
\]
and finally
\[\deg_z(M) \leq  r^2(\trs+1) \mathcal{E} + \trs r+ \frac12 r^2(r-1) (\trs-1) 
.\]
This concludes the proof of Inequality \eqref{eq:n16011} in Theorem~\ref{theo:1} in the Fuchsian case, i.e., when $\mathcal{N}=0$.

\subsection{General case}\label{ssec:casgeneral}

Again, we follow a strategy similar to that of van~Hoeij~\cite{Hoeij97b}, replacing Fuchs' relation by a generalization due to Bertrand and Laumon~\cite[Appendix]{Bertrand88}, see also~\cite{corel}.

\subsubsection*{Newton polygons} Part of the information on the degrees of factors comes from patching up local information at each singularity that can be read off the Newton polygons of the operators. We first recall their main definitions and properties (see \cite[p. 90, \S 3.3]{putsinger}).
Let
\begin{equation}\label{eq:5*}
R=\sum_{j=0}^\mu c_j(z) \partial_z^j=\sum_{j=0}^\mu\frac{U_j(z)}{V(z)}\partial_z^j \in \mathbb C(z)[\partial_z],
\end{equation}
with $U_j,V\in \mathbb C[z]$ and $V$ of lowest degree and $c_\mu$ not necessarily equal to~1. 
Consider the polynomials  $P_0(z)\equiv 1$ and, for $j\ge 1$,   $P_j(z):=\prod_{k=0}^{j-1}(z-k)\in \mathbb C[z]$ of degree $j$: they are such that $\partial_z^j = z^{-j} P_j(\theta_z)$ for all $j\ge 0$, where $\theta_z:=z\partial_z$. 
Rewriting $R=\sum_{j=0}^\mu c_j(z)z^{-j}P_j(\theta_z)$, the \emph{Newton polygon of $R$ at~$0$} is obtained by taking the lower-left boundary of the convex hull
of the points $(j,i)\in\mathbb{R}^2$ such that the coefficient of
$z^i\theta^j$ in the Laurent expansion of~$R$ at~$z=0$ is nonzero. The Newton polygon at another finite
point~$\rho\in\mathbb{C}$ is obtained similarly with
$\theta_{\rho,z}=(z-\rho)\partial_z$ and coefficients in $\mathbb{C}((z-\rho))$,
while the Newton polygon at infinity is the Newton polygon at~$0$ of the
operator $\widetilde{R}$ obtained from $R$ by changing $z$ into $1/z$. By
definition, the slopes of the Newton polygon  at~$\rho$ (finite or not) are all $\ge 0$ and
they are $0$ if and only if $R$ is regular or regular singular at~$\rho$.

In this work, we only use the largest slope of $R$ at~$\rho \in \mathbb
C\cup\{\infty\}$, that we denote by~$\mathcal{N}_{\rho}(R)$; this is also 
known as the {\em Katz rank} of $R$ at $\rho$, see~\cite[pp.~229--231]{BertrandBBKI} and~\cite{Bertrand88}. 
When $L=NM$,  
for any  $\rho \in \mathbb C\cup \{\infty\}$, we have 
\begin{equation} \label{ineq_slopes}
\mathcal{N}_{\rho}(M)\le \mathcal{N}_{\rho}(L).
\end{equation}
Indeed, a fundamental property is that the Newton polygon
of a product of operators is the Minkowski sum of their Newton
polygons~\cite[p. 92, Lemma~3.45]{putsinger}. Hence, the slopes of
$M$ at any point~$\rho$ form a subset of those of $L$ at~$\rho$. 

\bigskip

We now assume $R\in \mathbb C(z)[\partial_z]$ to be monic and of the form \eqref{eq:5*}. 
Let 
$v_j:=\textup{val}_{z=0}(c_j(z))$ for $j\leq\mu$. 
Note that $v_\mu=0$. 
Then for any $j\in \{0,\ldots, \mu-1\}$, we have  
\begin{equation}\label{eq:t00}
\mathcal{N}_{0}(R) \ge \frac{(v_\mu-\mu)-(v_j-j)}{\mu-j}=-1-\frac{v_j}{\mu-j}.
\end{equation}
It follows that for any $j\in \{0,\ldots, \mu-1\}$,
\begin{equation*}
\textup{val}_{z=0}\big(c_j(z)\big) \ge -\mu (\mathcal{N}_{0}(R) +1).
\end{equation*}
By a similar reasoning, for any finite $\rho \in \mathbb C$ and for any $j\in \{0,\ldots, \mu-1\}$,
\begin{equation}\label{eq:t2}
\textup{val}_{z=\rho}\big(c_j(z)\big) \ge -\mu (\mathcal{N}_{\rho}(R)+1).
\end{equation}
To deal with $\rho=\infty$, we remark that setting $\zeta=1/z$, we have $\zeta\partial_{\zeta}=-z\partial_z$. Hence, 
$$
\widetilde{R} = \sum_{j=0}^\mu \big(z^{j} c_j(1/z)\big)Q_j(\theta_z)
$$
where $Q_j(X):=P_j(-X)$. 
In view of $\operatorname{val}_{z=0} (z^{j}c_j(1/z))=j-\deg(c_j(z))$, the analogue of the inequality~\eqref{eq:t00} is then 
\[\mathcal{N}_\infty(R)\ge\frac{\mu-(j-\deg(c_j))}{\mu-j}=1+\frac{\deg(c_j)}{\mu-j},\qquad j=0,\dots,\mu-1,\]
leading to the bound
\begin{equation}\label{eq:t3}
\deg(c_j)  \le \mu (\mathcal{N}_{\infty}(R)-1), \qquad j=0, \ldots, \mu-1.
\end{equation}
Any finite singularity $\rho$ of $R$ is a  root of $V$ and there exists $j_\rho\in \{0, \ldots, \mu-1\}$ such that $\rho$ is not a root of $U_{j_\rho}$, so that  $\textup{val}_{z=\rho}\big(c_{j_\rho}(z)\big)=-\textup{val}_{z=\rho}\big(V(z)\big)$. Using \eqref{eq:t2} with $j=j_\rho$, we thus deduce that
\begin{equation}\label{eq:t4}
\textup{val}_{z=\rho}\big(V(z)\big) \le \mu (\mathcal{N}_{\rho}(R) +1).
\end{equation}
A similar reasoning at infinity gives
\begin{equation}\label{eq:t5}
\deg(U_j) \le \deg(V) + \mu (\mathcal{N}_{\infty}(R) -1), \qquad j=0, \ldots, \mu-1.
\end{equation}

\medskip

Let now $L=NM$ be a factorization of $L$ with a monic factor $M\in \mathbb C(z)[\frac{d}{dz}]$. We apply the bounds above to $R:=M$ and
$\mu:=r$.  Set
\[
\mathcal{N}=\max_{\rho \in \sing(L)\cup \{\infty\}} \mathcal{N}_{\rho}(L)\qquad\text{and} 
\qquad M=\sum_{j=0}^r \frac{A_j(z)}{B(z)}\partial_z^j,
\]
where the $A_j$'s and $B$ are as in \eqref{eq:5*}. In particular, by~\eqref{ineq_slopes}, 
for any $j=0, \ldots, r-1$, 
\begin{align}
\deg(A_j) &\le \deg(B) + r \mathcal{N}_{\infty}(M) -r \notag
\\
&\le \deg(B) + r \mathcal{N}-r. \label{eq:t6}
\end{align}
If $\rho \in \sing(M)$, then Eq.~\eqref{eq:t4} gives 
\begin{equation}\label{eq:t7}
\textup{val}_{z=\rho}\big(B(z)\big) \le r (\mathcal{N}_{\rho}(M) +1).
\end{equation}
If furthermore $\rho \in \alpha(M)$, then in particular it is a regular singularity, so that $\mathcal{N}_{\rho}(M)=0$ and this bound reduces to
\begin{equation} \label{eq:t8}
\textup{val}_{z=\rho}\big(B(z)\big) \le r.
\end{equation}
Since the degree of a polynomial is the sum of the valuations (or multiplicities) at its roots, it follows from \eqref{eq:t7}, \eqref{eq:t8} and $\sigma(M)\subset \sigma(L)$ that 
\begin{align}
\deg(B) &\le r( \mathcal{N} +1)\cdot\#\sigma(M)+
r\cdot \#\alpha(M) \notag
\\ 
&\le r( \mathcal{N} +1)\trs+ r \cdot \#\alpha(M). \label{eq:t9}
\end{align}
From \eqref{eq:t6} and  \eqref{eq:t9}, we see that an explicit upper bound for 
\begin{align}
\deg_z(M)&:=
\max(\deg(A_0), \ldots, \deg(A_{r-1}), \deg(B)) \notag 
\\
&\le \deg(B) + r \mathcal{N} \label{eq:t10}
\\
&\le r\trs + r(\trs+1) \mathcal{N} + r\cdot \#\alpha(M) \label{eq:t100}
\end{align}
will again be obtained from an explicit upper bound for $\#\alpha(M)$. (If  $\mathcal{N}\ge 1$, then the right-hand side of \eqref{eq:t10} can be improved to  $\deg(B) + r (\mathcal{N} -1)$ by \eqref{eq:t6}, with corresponding improvements in subsequent inequalities.)

\subsubsection*{Generalized Fuchs' relation} For any $R\in \mathbb C(z)[\partial_z]$ of order $r$, we consider the $D$-module $\widehat{R}:=\mathbb C(z)[\partial_z]/(\mathbb C(z)[\partial_z]R)$. 
The generalization of Fuchs' relation~\eqref{eq:Fuchs} given in~\cite[Appendice, p.~84]{Bertrand88}, \cite[p. 53, Theorem 2]{Bertrand99} and \cite[p. 298]{corel} is
\begin{equation}\label{eq:Fuchs3}
\sum_{\rho\in \sing(R)\cup \{\infty\}} \Big(S_\rho(R)- \frac{1}{2}\textup{irr}_\rho\big(\textup{End}(\widehat{R})\big)\Big)= -r(r-1),
\end{equation}
where as before 
\begin{equation}\label{eq:fuchs4}
S_\rho({R}) := \sum_{j=1}^{r} e_j(\rho) - \frac{r(r-1)}{2},
\end{equation}
but now the $e_j(\rho)$'s are the \emph{generalized local exponents} of $R$ at the point~$\rho \in \mathbb C \cup \{\infty\}$ (see \cite[Appendice, pp. 82-83]{Bertrand88} or \cite[p. 297]{corel} for their definition). 

Given a differential operator $R$ in $\mathbb{C}(z)[\partial_z]$, its {\em
Malgrange's irregularity}~\cite{Malgrange71}, denoted $\textup{irr}_\rho (\widehat{R})$,
is a non-negative integer which measures the defect of Fuchsianity of $R$ at~$\rho$. 
Next, $\textup{End}(\widehat{R})$ in \eqref{eq:Fuchs3} is isomorphic to the $D$-module $\widehat{R}\otimes \widehat{R^*}$, where $R^*$ is the adjoint of~$R$. 
By~\cite[Appendice, p. 84]{Bertrand88}, the integer
 $\textup{irr}_\rho\big(\textup{End}(\widehat{R})\big)$ can be bounded in terms
of $\mathcal{N}_{\rho}(R)$: for any $\rho \in \mathbb C \cup \{\infty\}$, 
\begin{equation}\label{ineq_Bertrand88}
\textup{irr}_\rho\big(\textup{End}(\widehat{R})\big)\le r(r-1)\mathcal{N}_{\rho}(R).
\end{equation}
If $\rho\in \mathbb C \cup \{\infty\}$ is an ordinary point or a regular
singularity of $R$, we have $\mathcal{N}_{\rho}(R)=0$ and {\em a fortiori}  $\textup{irr}_\rho\big(\textup{End}(\widehat{R})\big)=0$ as well; 
we thus recover the usual Fuchs relation \eqref{eq:Fuchs} when $R$ is Fuchsian. 

We are now ready to bound $\#\alpha(M)$ in any factorization $L=NM$. We recall that $\alpha(M)$ and $\sigma(M)$ form a partition of $\sing(M)$, and that $\sigma(M)\subset \sigma(L)$.  
Therefore, 
$$
r(r-1) + 
\sum_{\rho \in \alpha(M)} S_{\rho}(M) = -\sum_{\rho \in \sigma(M) \cup \{\infty\}}S_{\rho}(M) + \frac{1}{2} \sum_{\rho \in \sigma(M)\cup \{\infty\}} \textup{irr}_\rho\big(\textup{End}(\widehat{M})\big).
$$
Now, $S_{\rho}(M)\in \mathbb N \setminus \{0\}$ for any $\rho \in \alpha(M)$ (again 
by \cite[Chap.~V, \S18, 
p.~69]{Poole1960}) 
and $\vert S_{\rho}(M)\vert \le r\mathcal{E}+\frac{r(r-1)}{2}$  for any $\rho \in \sigma(M)\cup \{\infty\}$ by \eqref{eq:fuchs4}.

It follows from \eqref{eq:newreferee} that 
\begin{align}
\#\alpha(M)&+r(r-1) \notag
\\
&\le (\#\sigma(M)+1)\Big(r\mathcal{E}+\frac{r(r-1)}{2}\Big)+
\frac{r(r-1)}{2} \sum_{\rho \in \sigma(M)\cup \{\infty\}} \mathcal{N}_{\rho}(M) \notag
\\
&\le (\#\sigma(L)
+1)\Big(r\mathcal{E}+\frac{r(r-1)}{2}\Big)+
\frac{r(r-1)}{2} \sum_{\rho \in \sigma(L)\cup \{\infty\}} \mathcal{N}_{\rho}(L) \notag
\\
&\le (\trs+1)\Big(r\mathcal{E}+\frac{r(r-1)}{2}\Big)+ \frac12r(r-1)(\trs+1)\mathcal{N}. 
\label{eq:t1506}
\end{align}
In this sequence of inequalities, the first one follows from~\eqref{ineq_Bertrand88} and 
the second one uses~\eqref{ineq_slopes}. 
Hence, 
\begin{equation*}
\#\alpha(M) \le  (\trs+1)r\mathcal{E}+ \frac 12r(r-1)(\trs+1)(\mathcal{N}+1)-r(r-1).
\end{equation*}
It follows  from \eqref{eq:t100} that 
\begin{equation}\label{eq:t11}
\qquad \deg_z(M)\le r^2(\trs+1)\mathcal{E}+ r(\mathcal{N}+1)\trs + r\mathcal{N} +\frac12 r^2(r-1)\big((\trs+1)(\mathcal{N}+1)-2\big). \qquad
\end{equation} 
This completes the proof of Theorem \ref{theo:1}.

\medskip

We have mentioned above that if $\mathcal{N}\ge 1$, then the term $r\mathcal{N}$ in \eqref{eq:t11} can be replaced by $r(\mathcal{N}-1)$. It is possible to improve further this bound. Indeed, for any $R\in \mathbb C(z)[\partial_z]$, we have 
\begin{equation*} 
\sum_{\rho
\in\sing(R) \cup \{\infty\} } \big(\mathcal{N}_{\rho}(R)+1\big) \le
2\deg_z(R)+2, 
\end{equation*}
by the arguments used in the proof of \cite[p. 185, Lemme~2bis]{BeBe85}. Hence, the final term $\frac12r(r-1)(\trs+1)\mathcal{N}$ in 
\eqref{eq:t1506} could be replaced by 
$$
\frac{r(r-1)}{2}\min\Big((\trs+1)\mathcal{N}, 2q+1-\#\textup{Sing}(L)\Big),
$$
with a corresponding improvement of \eqref{eq:t11}.
\medskip

It may seem at first sight that it should  be possible to adapt the proof of Theorem \ref{theo:1} to any algebraically closed field~$\mathbb K$ of characteristic 0, instead of $\mathbb C$. This is the case as long as {\em equalities} are used. However, the deductions made from the generalized Fuchs relation (which holds for such a field $\mathbb K$) are based on various {\em inequalities}. The argument might in principle be adapted when $\mathbb K$ is also endowed with an archimedean absolute value. But by Ostrowski's Theorem \cite[p. 33, Theorem 1.1]{cassels}, such a field can be embedded into a subfield of $\mathbb C$ endowed with an absolute value given by a positive power of the modulus, and we would in fact gain nothing. Finally, it is not clear to us how a bound similar to~\eqref{eq:t11} could be obtained with this method for a field $\mathbb K$ not endowed with an archimedean absolute value.
We point to~\cite[Sec.~9]{Hoeij97b} for a possible reduction of the general case to the one treated here.

\section*{Acknowledgements} We warmly thank Daniel Bertrand for his comments on a previous version of this article, and Mark van Hoeij for a question that led us to strengthen our main result. 
Bruno Salvy has been supported in part by FastRelax ANR-14-CE25-0018-01.

\end{document}